\documentclass[12pt, reqno]{amsart}
\usepackage{amsmath, amsthm, amscd, amsfonts, amssymb, graphicx, xcolor}
\usepackage[bookmarksnumbered, colorlinks, plainpages]{hyperref}

\newtheorem{theorem}{Theorem}[section]

\theoremstyle{definition}

\newtheorem{example}[theorem]{Example}

\theoremstyle{remark}

\numberwithin{equation}{section}

\begin{document}
\setcounter{page}{1}

\color{darkgray}{
\noindent 

\centerline{}

\centerline{}

\title[Sums and differences of sets ( improvement over AlphaEvolve )]{Sums and differences of sets ( improvement over AlphaEvolve )}

\author[Robert Gerbicz]{Robert Gerbicz}

\address{$^{1}$ Halasztelek, Hungary.}
\email{\textcolor[rgb]{0.00,0.00,0.84}{robert.gerbicz@gmail.com}}




\date{May 22, 2025.}

\begin{abstract}
On May 14, 2025, DeepMind announced that AlphaEvolve, a large language model applied to a set of mathematical problems, had matched or exceeded the best known bounds on several problems. In the case of the sum and difference of sets problem, AlphaEvolve, using a set of $54265$ integers, improved the known lower bound of $\theta=1.14465$ to $\theta=1.1584$. In this paper, we present an improved bound $\theta=1.173050$ using an explicit construction of a U set that contains more than $10^{43546}$ elements. For fast integer and floating-point arithmetic, we used the (free) GMP library.
\end{abstract} \maketitle

\section{Introduction and preliminaries}

\noindent
Let $K>1$ be a real number. In \cite{Gy2007} the authors asked for the largest possible value of $\theta$ such that there exist integer sets $A$ and $B$, with $|A|$ arbitrarily large, satisying the condition that for some constant $c(K)>0$:
\begin{equation}\label{1}
|A+B|\leq K|A|~\text{and}~|A-B|\geq c(K)\cdot |A+B|^{\theta}
\end{equation}

They showed that $\theta\geq 1+\frac{\log{(\frac{|U-U|}{|U+U|})}}{\log{(2\cdot \max(U)+1)}}$ where $U$ is a set of non-negative integers containing zero.

In their constructions, they started from the set $V=V(m,L)=\{(x_{1},\ldots,x_{m})\in \mathbb{N}^m:x_{1}+\ldots +x_{m}\leq L\}$, where $m\geq 0$ and $L\geq 0$. From combinatorics, it is well known that $|V(m,L)|=\binom{m+L}{m}$.

To obtain the integers of the set, they used the function $f:(x_{1},\ldots x_{m})\rightarrow \sum_{k=0}^{m-1}x_{k}\cdot L_{k}$, where $L_{0}=1$ and $L_{k}=2L\cdot L_{k-1}+1$ for $k>0$. It is easy to see that the function $f$ is injective on both $V+V$ and $V-V$. The required integers of the $U$ set were then defined by $U=\{f(x):x\in V(m,L)\}$.

Using formulas for $|V-V|,~|V+V|$ and $\max(U)$ they were able to establish a lower bound of $\theta=1.14465$. The known upper bound is $\theta\leq\frac{4}{3}$. DeepMind's AlphaEvolve (see \cite{N2025}), using a large language model, improved the lower bound to  $\theta=1.1584$. 

My idea was based on the observation that for $x_{1}+\ldots +x_{m}=s\leq L$, the average value of each $x_{i}$ is only $\frac{s}{m}\leq \frac{L}{m}$. So for vectors $x$ and $x'$, the average value of $x_{k}+x'_{k}$ is still relatively small. This suggests that we might improve the value of $\theta$ by reducing the large $L_{k}$ values, which are roughly $(2L)^{k}$. To achieve this, we restrict the coordinates by requiring $x_{k}\leq B$, and then we can use $L_{k}=(2\cdot B+1)^k$. This ensures that the function $f$ remains injective on the sets $V-V$ and $V+V$.

\section{Main results}

Let $B\geq 0$ be an integer. Define:
\begin{equation}\label{2}
W=W(m,L,B)=\{x=(x_{1},\ldots,x_{m})\in\mathbb{N}^m\land x_{1}\leq B,\ldots ,x_{m}\leq B\}
\end{equation}
Clearly, if $m=0\lor L=0\lor B=0$, then $|W(m,L,B)|=1$. Therefore, in the following formulas, we may assume $m>0,~L>0,~B>0$. For $1\leq k\leq m$, let $H_{k}=\{x:x=(x_{1},\ldots,x_{m})\in W(m,L,B)\land x_{k}\geq B+1\}$. We easily have $|H_{1}\cap \ldots \cap H_{k}|=|V(m,L-k(B+1))|=\binom{m+L-k(B+1)}{m}$, since for any $x\in H_{1}\cap\ldots\cap H_{k}$, the first $k$ coordinates are each at least $B+1$. Thus, we are left with a coordinate sum of at most $L-k(B+1)$ for the shifted vector $(x_{1}-(B+1),\ldots,x_{k}-(B+1),x_{k+1},\ldots,x_{m})$. Because this holds for any choice of $k$ coordinate indices, using the inclusion-exclusion principle, we obtain:
\begin{align}\label{3}
|W(m,L,B)|=|V(m,L)\setminus (H_{1}\cup \ldots \cup H_{m})|=\nonumber\\
|V(m,L)|+\sum_{k=1}^{\lfloor\frac{L}{B+1}\rfloor}(-1)^{k} \binom{m}{k}\binom{m+L-k(B+1)}{m}=\nonumber\\
\sum_{k=0}^{\lfloor\frac{L}{B+1}\rfloor}(-1)^{k} \binom{m}{k}\binom{m+L-k(B+1)}{m}\end{align}
So we have:
\begin{equation}\label{4}
|W(m,L,B)|=\sum_{k=0}^{\lfloor\frac{L}{B+1}\rfloor}(-1)^{k} \binom{m}{k}\binom{m+L-k(B+1)}{m}\end{equation}

To obtain the required set of integers $U$, consider the function\newline
$g:(x_{1},\ldots,x_{m})\rightarrow \sum_{k=0}^{m-1}x_{k}(2\cdot B+1)^{k}$, and define $U=U(m,L,B)=\{g(x):x\in W(m,L,B)\}$. It is clear that $g$ is injective on both $W+W$ and $W-W$.

We have:
\begin{equation}\label{5}
s(U):=|U+U|=|W+W|=|W(k,2L,2B)|
\end{equation}
since for each vector in $W+W$, the coordinate sum is at most $2L$, and each coordinate is at most $2B$.

To compute $|W-W|$, observe that $x\in W-W$ if and only if the absolute value of each coordinate is at most $B$, the sum of  positive coordinates is at most $L$, and the absolute value of the sum of the remaining non-positive coordinates is also at most $L$. Suppose there are $k$ positive coordinates. Then, we can choose their positions in $\binom{m}{k}$ ways. There are $|W(k,L-k,B-1)|$ ways to set these $k$ coordinates. For the remaining $m-k$ coordinates, there are $|W(m-k,L,B)|$ choices to assign values to the non-positive coordinates. Since these choices are independent:
\begin{equation}\label{6}
d(U):=|U-U|=|W-W|=\sum_{k=0}^{\min(m,L)}\binom{m}{k}\cdot |W(k,L-k,B-1)|\cdot |W(m-k,L,B)|
\end{equation}

Since each coordinate is at most $B$, we may assume $L\leq mB$. The maximum of $U$, denoted $\max(U)$, is obtained when the coordinate sum reaches its maximum, namely $L$. So in a greedy way set the largest coordinates to $B$, while we can, then set a coordinate to $L\%B$, and set the rest to zero.

For example, if $m=4,~L=8,~B=3$, then $\max(U)=3\cdot 7^{3}+3\cdot 7^{2}+2\cdot 7^{1}=1190$. In general, using  $t=\lfloor\frac{L}{B}\rfloor$, we can compute:
\begin{equation}\label{7}
q(U):=2\max(U)+1=(2B+1)^{m}-(2B+1)^{m-t}+2(L\%B)(2B+1)^{m-t-1}+1
\end{equation}

\begin{example}
For $m=4,~ L=8,~ B=3$, we have $U=U(4,8,3)$ with: $|U|=221,s(U)=|U+U|=2075,d(U)=|U-U|=2307,q(U)=2\max(U)+1=2381$. These values yield: $\theta\geq 1+\frac{\log{(\frac{|U-U|}{|U+U|})}}{\log{(2\cdot \max(U)+1)}}=1.013631$.
\end{example}

\begin{theorem}
For the sums and differences of sets problem, we have\newline $\theta\geq 1.173050$.
\end{theorem}
\textbf{Proof.} Larger values of $m$ and $L$ yield very large integers for $s(U),~d(U)$ and $|U|$. Therefore, we first use floating-point arithmetic to search for good candidates yielding larger values of $\theta$, and then confirm these values using integer arithmetic.

The formula in \eqref{4} can result in large numerical errors when using floating-point arithmetic, as it involves an alternating sum. Grouping positive and negative terms does not sufficiently mitigate this precision issue. Thus, when only an approximation of $W(m,L,B)$ is needed, we compute it exactly using integer arithmetic and then round the result.

To find a better value than AlphaEvolve's $\theta=1.1584$, we performed a somewhat restricted search with $1\leq m\leq 128,~ L=64$ and $1\leq B\leq 7$. In less than one second (on a single CPU core), the search revealed that the maximum value, $\theta=1.162997$, is achieved  for $m=80,~L=64,~B=5$. Integer arithmetic confirms this value of $\theta$.

In these small-scale searches, we observed that for fixed $L$, the largest $\theta$ value tends to occur when $B=5$ and $m\approx \frac{5}{4}L$. In the above example, this relation holds exactly, though it is not universally true. For $L=2^{16}=65536$, we searched values of $m$ near $\frac{5}{4}\cdot 65536=81920$, and found that (via floating point estimation) that $m=81411$ likely provides the best $\theta$ value. Using exact integer arithmetic, we confirmed in approximately $15$ hours that $\theta=1.173050$ for $m=81411,~L=65536,~B=5$. Rounding the calculated exact values gives:
\begin{align*}
|U|\approx 6.314107319\cdot 10^{43546}\\
s(U)=|U+U|\approx 3.208492702\cdot 10^{61228}\\
d(U)=|U-U|\approx 6.587554451\cdot 10^{75899}\\
q(U)=2\max(U)+1\approx 6.605282799\cdot 10^{84780}
\end{align*}

For efficient integer and floating-point arithmetic, we used the (free) GMP library: \cite{GMP}. The GMP code and the exact computed values for the triple $m=81411,~L=65536,~B=5$ are available in \cite{Gerbicz}.

Further improvements might be possible by increasing the value of $L$, but any additional gain is likely to be less than $0.0001$ beyond the achieved value of $\theta=1.173050$.

\bibliographystyle{amsplain}

\begin{thebibliography}{99}

\bibitem{Gerbicz} Sums and differences of sets. \url{https://bit.ly/43mWCVC}

\bibitem{GMP} The GNU Multiple Precision Arithmetic Library. \url{https://gmplib.org/}

\bibitem{Gy2007} K. Gyarmati, F. Hennecart, and I. Z. Ruzsa (2007). Sums and differences of finite sets. Functiones et Approximatio Commentarii Mathematici, 37(1):175–186, 2007. \url{https://gyarmatikati.web.elte.hu/publ/sumdiffv.pdf}

\bibitem{N2025}  Alexander Novikov, Ngân Vu, Marvin Eisenberger, Emilien Dupont, Po-Sen Huang, Adam Zsolt Wagner, Sergey Shirobokov, Borislav Kozlovskii, Francisco J. R. Ruiz, Abbas Mehrabian, M. Pawan Kumar, Abigail
See, Swarat Chaudhuri, George Holland, Alex Davies, Sebastian Nowozin, Pushmeet Kohli and Matej Balog (2025). AlphaEvolve: A coding agent for scientific and algorithmic discovery. 
\url{http://bit.ly/4jkeOoX}
\end{thebibliography}

\end{document}